\documentclass{conm-p-l}
\pdfoutput=1
\usepackage{amscd}
\usepackage{amssymb}
\usepackage[all]{xy}
\usepackage{graphicx}
\usepackage[T1]{fontenc}
\usepackage{lmodern}   
\date{26 July 2011}
\title[Twisted Deformation Quantization]{Twisted Deformation
Quantization of Algebraic Varieties (Survey)}

\author{Amnon Yekutieli}
\address{Department of  Mathematics,
Ben Gurion University,
Be'er Sheva 84105,
Israel}
\email{amyekut@math.bgu.ac.il}
\thanks{
Primary: 53D55; Secondary: 14B10, 16S80, 17B40, 18D05.}
\thanks{This research was supported by the US-Israel Binational
Science Foundation and by the Israel Science Foundation.}

\newtheorem{thm}[equation]{Theorem}
\newtheorem{cor}[equation]{Corollary}

\theoremstyle{definition}
\newtheorem{dfn}[equation]{Definition}
\newtheorem{rem}[equation]{Remark}
\newtheorem{exa}[equation]{Example}

\newtheorem{que}[equation]{Question}

\numberwithin{equation}{section}

\newcommand{\ol}{\overline}
\newcommand{\iso}{\stackrel{\simeq}{\longrightarrow}}

\newcommand{\opn}{\operatorname}
\newcommand{\cat}[1]{\operatorname{\mathsf{#1}}}

\newcommand{\rmitem}[1]{\item[\text{\textup{(#1)}}]}
\newcommand{\mfrak}[1]{\mathfrak{#1}}
\newcommand{\mcal}[1]{\mathcal{#1}}

\newcommand{\mbf}[1]{\mathbf{#1}}
\newcommand{\mrm}[1]{\mathrm{#1}}
\newcommand{\mbb}[1]{\mathbb{#1}}

\newcommand{\tup}[1]{\textup{#1}}
\newcommand{\bsym}[1]{\boldsymbol{#1}}
\newcommand{\bwedge}{\bigwedge\nolimits}
\newcommand{\boplus}{\bigoplus\nolimits}

\newcommand{\what}{\widehat}

\newcommand{\til}[1]{\tilde{#1}}
\renewcommand{\d}{\mathrm{d}}

\DeclareMathSymbol{\mathbbk}{\mathord}{AMSb}{"7C}

\newcommand{\K}{\mbb{K}}
\newcommand{\R}{\mbb{R}}
\newcommand{\m}{\mfrak{m}}
\newcommand{\g}{\mfrak{g}}
\newcommand{\h}{\mfrak{h}}
\newcommand{\smfrac}[2]{\textstyle \frac{#1}{#2}}
\newcommand{\hatotimes}[1]{\, \what{\otimes}_{#1} \,}
\newcommand{\gerbe}[1]{\bsym{\mcal{#1}}}

\begin{document}

\begin{abstract}
Let $X$ be a smooth algebraic variety over a field of characteristic
$0$. We introduce the notion of twisted associative (resp.\
Poisson) deformation of the structure sheaf $\mcal{O}_X$. These are
stack-like versions of usual deformations. We prove that there is a
twisted quantization map from twisted Poisson deformations
to twisted associative deformations, which is canonical and bijective
on equivalence classes.
\end{abstract}

\maketitle

\setcounter{section}{-1}
\section{Introduction}
\label{sec:Int}

Let $\K$ be a field of characteristic $0$, and let $X$ be a
smooth algebraic variety over $\K$, with structure sheaf
$\mcal{O}_X$. An {\em associative deformation of $\mcal{O}_X$} is a
sheaf $\mcal{A}$ of flat complete associative $\K[[\hbar]]$-algebras on $X$,
with an isomorphism $\K \otimes_{\K[[\hbar]]} \mcal{A} \cong \mcal{O}_X$, called
an augmentation. Similarly, a {\em Poisson deformation of $\mcal{O}_X$} is a
sheaf $\mcal{A}$ of flat complete commutative Poisson $\K[[\hbar]]$-algebras
on $X$, with an augmentation to $\mcal{O}_X$. 

In this paper we introduce the notion of {\em twisted associative
\tup{(}resp.\ Poisson\tup{)} deformation of $\mcal{O}_X$}.  
A twisted deformation (or either kind) is a stack-like version of an
ordinary deformation. The precise definition is given in Section 
\ref{sec:twdefs}. But to give an idea, let us say that a twisted deformation 
$\gerbe{A}$ can be described as a collection of locally defined
deformations $\mcal{A}_i$, each living on an open set $U_i$ of $X$, that
are glued together in a loose way. We should also say that an
associative deformation $\gerbe{A}$ is a special kind of
{\em stack of algebroids}, in the sense of \cite{Ko2}.
Indeed, one reason for introducing twisted deformations is to have a
Poisson analogue of a stack of algebroids.

There is a notion of {\em twisted gauge equivalence} between  
twisted associative (resp.\ Poisson) deformation of
$\mcal{O}_X$. A twisted deformation $\gerbe{A}$ induces a {\em
first order bracket} $\{ -,- \}_{\gerbe{A}}$ on $\mcal{O}_X$.
The main result is Theorem \ref{thm_new}, which says that 
there is a canonical bijection of sets
\[ \begin{aligned}
& \opn{tw{.}quant} : 
\frac{ \{ \tup{twisted Poisson deformations of 
$\mcal{O}_X$} \} }
{\tup{ twisted gauge equivalence}} \\
& \qquad \qquad \iso \
\frac{ \{ \tup{twisted associative deformations of 
$\mcal{O}_X$} \} }
{\tup{ twisted gauge equivalence}}
\end{aligned} \] 
called the {\em twisted quantization map}. 
It preserves first order brackets, and commutes with \'etale morphisms 
$X' \to X$.

This survey article is an edited version of lectures that I gave on several
occasions. In addition to the main body of the article, there are four
appendices, that provide further details on certain aspects of the topic. 
Full details can be found in my papers listed in the bibliography. The
bibliography also covers work by other researchers in this area.

\medskip \noindent
\textbf{Acknowledgments.}
Part of the work presented here is joint with Frederick Leitner.
Many of the ideas in this paper are influenced by the work
of Maxim Kontsevich, and I wish to thank him for discussing this
material with me.
Thanks also to Michael Artin, Pavel Etingof, Damien
Calaque, Michel Van den Bergh, Pierre Deligne, 
Lawrence Breen, Pierre Schapira and James Stasheff 
for their assistance on various aspects of the paper.

\section{Some background on Deformation Quantization}

Throughout $\K$ is a field of characteristic $0$.

Let $C$ be a commutative $\K$-algebra. 
Recall that a {\em Poisson bracket} on $C$ is a 
$\K$-bilinear function
\[ \{-,-\} : C \times C \to C \]
which makes $C$ into a Lie algebra, and is a biderivation (i.e.\ 
a derivation in each argument). 
The pair $\bigl( C, \{-,-\} \bigr)$ is called a Poisson algebra. 

Poisson algebras arise in several ways, e.g.\ classical Hamiltonian
mechanics, or Lie theory.

Let $\K[[\hbar]]$ be the ring of formal power series in the
variable $\hbar$. And let $C[[\hbar]]$ be the set of formal power series with
coefficients in $C$, which we view only as a 
$\K[[\hbar]]$-module. 
A {\em star product} on $C[[\hbar]]$ is a $\K[[\hbar]]$-bilinear function
\[ \star : C[[\hbar]] \times C[[\hbar]] \to C[[\hbar]]  \]
which makes $C[[\hbar]]$ into an associative 
$\K[[\hbar]]$-algebra, with unit $1 \in C$, and  such that 
\[ f \star g \equiv f g \ \opn{mod}\ \hbar \]
for any $g, f \in C$.

Note that the star product $\star$ can be expanded into a power series as
follows: there is a sequence $\{ \omega_j \}_{j \geq 1}$ of $\K$-bilinear
functions $\omega_j : C \times C \to C$, such that 
\begin{equation} \label{eqn:omega2}
f \star g := f g + 
\sum_{j = 1}^{\infty} \omega_j(f, g) \hbar^j
\end{equation}
for $f, g \in C$. The conditions that $\star$ is associative and unital place
certain constraints on the sequence $\{ \omega_j \}_{j \geq 1}$.

\begin{exa} \label{exa3}
Suppose $\star$ is a star product on  $C[[\hbar]]$. Given $f, g \in C$, we know
that 
\[ f \star g - g \star f \equiv 0 \ \opn{mod}\ \hbar . \]
Hence there is a unique element 
$\{ f, g \}_{\star} \in C$
such that 
\[ \smfrac{1}{2 \hbar} \bigl( f \star g - g \star f \bigr)
\equiv \{ f, g \}_{\star}  \ \opn{mod}\ \hbar . \]
It is quite easy to show that $\{ -,- \}_{\star}$ is a Poisson
bracket on $C$. We call it the {\em first order bracket}
of $\star$. 
\end{exa}

Deformation quantization seeks to reverse Example \ref{exa3}.

\begin{dfn}
Given a Poisson bracket $\{ -,- \}$ on the algebra $C$, a 
{\em deformation quantization} of $\{ -,- \}$ is a star product $\star$ on
$C[[\hbar]]$ whose first order bracket is $\{ -,- \}$.
\end{dfn}

In physics $\hbar$ is the {\em Planck constant}. For a
quantum
phenomenon depending on $\hbar$, the limit as $\hbar \to 0$ is thought
of the as the classical limit of this phenomenon.

The original idea by the physicists Flato et.\ al.\ (\cite{BFFLS},
1978) was that deformation quantization should model the transition
from classical Hamiltonian mechanics to quantum mechanics.
Special cases (like the Moyal product) were known. The problem 
arose: {\em does any Poisson bracket admit a 
deformation quantization}?

For a symplectic manifold $X$ and $C = \mrm{C}^{\infty}(X)$ 
the problem was solved by De Wilde and Lecomte  (\cite{DL}, 1983). 
A more geometric solution was discovered by Fedosov (\cite{Fe},
1994). The general case, i.e.\ $C = \mrm{C}^{\infty}(X)$ for a 
Poisson manifold $X$, was solved by Kontsevich  (\cite{Ko1}, 1997).
See surveys in the book \cite{CKTB}.

\begin{rem} \label{rem:parameter}
The deformations considered in this paper are parameterized by the algebra 
$\K[[\hbar]]$. It is possible to replace $\K[[\hbar]]$ with any other noetherian
complete local commutative $\K$-algebra $R$, with maximal ideal $\m$, such that
$R / \m = \K$. Instead of $C[[\hbar]]$ we take the complete tensor product 
$R \hatotimes{\K} C$, with its obvious augmentation to $C$.
All results will remain valid. 
\end{rem}

\section{Poisson Deformations of Algebraic Varieties}

In algebraic geometry we have to consider deformations as sheaves. 

Let $X$ be a smooth algebraic variety over $\K$, with structure 
sheaf $\mcal{O}_{X}$. 
We view $\mcal{O}_{X}$ as a Poisson $\K$-algebra with zero bracket.

\begin{dfn}
A {\em Poisson deformation} of $\mcal{O}_X$ 
is a sheaf $\mcal{A}$ of flat, $\hbar$-adically complete,
commutative Poisson $\K[[\hbar]]$-algebras on $X$, with an 
isomorphism of Poisson algebras 
\[ \psi : \mcal{A} / (\hbar) \iso \mcal{O}_X , \]
called an augmentation. 

A {\em gauge equivalence} $\mcal{A} \to \mcal{A}'$
between Poisson deformations is a $\K[[\hbar]]$-linear 
isomorphism of sheaves of Poisson algebras, that commutes with the 
augmentations to $\mcal{O}_X$. 
\end{dfn}

It may happen that $\mcal{A} \cong \mcal{O}_X[[\hbar]]$ as sheaves of
commutative $\K[[\hbar]]$-algebras augmented to $\mcal{O}_X$; if this is so,
then we say that $\mcal{A}$ is a sheaf-theoretically trivial deformation.
A sufficient condition for that is the vanishing of the cohomology group 
$\mrm{H}^1(X, \mcal{T}_X)$, where $\mcal{T}_X$ is the tangent sheaf. 
The corresponding $\K[[\hbar]]$-bilinear Poisson bracket on
$\mcal{O}_X[[\hbar]]$ is called a {\em formal Poisson bracket}.

Given a Poisson deformation $\mcal{A}$ of $\mcal{O}_X$, we 
may define the {\em first order bracket}
\[ \{ -,- \}_{\mcal{A}} :  \mcal{O}_X \times \mcal{O}_X \to \mcal{O}_X
. \]
This is a Poisson bracket whose formula is
\[ \{ f,g \}_{\mcal{A}} := \psi \bigl( \smfrac{1}{\hbar} \{ \til{f},
\til{g} 
\} \bigr) , \]
where $f, g \in \mcal{O}_X$ are local sections, and 
$\til{f}, \til{g} \in \mcal{A}$ are arbitrary local lifts. %
The first order bracket  is invariant under gauge equivalences. 

\begin{exa} \label{exa10}
Let $\{ -,- \}_1$ be some Poisson bracket on $\mcal{O}_{X}$. 
Put on the commutative $\K[[\hbar]]$-algebra 
$\mcal{A} := \mcal{O}_X [[\hbar]]$
the formal Poisson bracket
$\hbar \{ -,- \}_1$, namely 
\[ \{ f, g \} = \hbar \{ f, g \}_1  \]
for $f, g \in \mcal{O}_X$. 
Then $\mcal{A}$ is a Poisson deformation of $\mcal{O}_X$. The first
order bracket in this case is just
\[ \{ -,- \}_{\mcal{A}} = \{ -,- \}_{1} . \]
\end{exa}

Poisson deformations are controlled by a coherent sheaf of DG (differential
graded) Lie algebras 
$\mcal{T}^{}_{\mrm{poly}, X}$,
called the {\em poly derivations}. 
This is explained in Appendix A.

\section{Associative Deformations of Algebraic Varieties}

Let $X$ be a smooth algebraic variety over $\K$ as before.

\begin{dfn}
An {\em associative deformation of $\mcal{O}_X$} 
is a sheaf $\mcal{A}$ 
of flat, $\hbar$-adically complete, associative, unital 
$\K[[\hbar]]$-algebras on $X$, with an isomorphism of algebras
\[ \psi : \mcal{A} / (\hbar) \iso \mcal{O}_X , \]
called an augmentation.

A {\em gauge equivalence} $\mcal{A} \to \mcal{A}'$
between associative deformations is a $\K[[\hbar]]$-linear 
isomorphism of sheaves of unital algebras, that commutes with the 
augmentations to $\mcal{O}_X$. 
\end{dfn}

It may happen that $\mcal{A} \cong \mcal{O}_X[[\hbar]]$ as sheaves of
$\K[[\hbar]]$-modules augmented to $\mcal{O}_X$; if this is so,
then we say that $\mcal{A}$ is a sheaf-theoretically trivial deformation.
A sufficient condition for that is the vanishing of the cohomology group 
$\mrm{H}^1(X, \mcal{D}_X)$, where $\mcal{D}_X$ is the sheaf of differential
operators on $X$. (This fact is quite hard to prove, and it relies on Theorem
\ref{thm:diff}.)
The corresponding $\K[[\hbar]]$-bilinear multiplication on
$\mcal{O}_X[[\hbar]]$ is called a {\em star product}, like in Section 1.

Given an associative deformation $\mcal{A}$ we 
may define the first order bracket
\[ \{ -,- \}_{\mcal{A}} :  \mcal{O}_X \times \mcal{O}_X \to \mcal{O}_X
. \]
The formula, in terms of local sections, is
\[ \{ f,g \}_{\mcal{A}} :=
\psi \bigl( \smfrac{1}{2 \hbar} (\til{f} \star 
\til{g} - \til{g} \star \til{f}) \bigr) . \]
The first order bracket is invariant under gauge equivalences.

Associative deformations are controlled by a quasi-coherent sheaf of
DG Lie algebras $\mcal{D}^{}_{\mrm{poly}, X}$,
called the {\em poly differential operators}. 
This is explained in Appendix A.

Note that both kinds of deformations -- Poisson and associative --
include as special cases the classical commutative deformations of $\mcal{O}_X$.

\section{Deformation Quantization}

Kontsevich \cite{Ko1} proved that any Poisson deformation of a real 
$\mrm{C}^{\infty}$ manifold $X$ can be canonically quantized. 
In this section we present an algebraic version of this result. But first a
definition. 

\begin{dfn}
Let $\mcal{A}$ be a Poisson deformation of $\mcal{O}_X$. %
A {\em quantization} of $\mcal{A}$ is an 
associative deformation $\mcal{B}$, such that the first order
brackets satisfy
\[ \{ -,- \}_{\mcal{B}} = \{ -,- \}_{\mcal{A}} . \]
\end{dfn} 

Recalling Example 2.2,
we see that this definition captures the essence
of deformation quantization, namely quantizing a Poisson bracket on 
$\mcal{O}_X$.

\begin{thm}[\cite{Ye1}] \label{thm1} 
Let $\K$ be a field containing $\mbb{R}$, and let $X$ be a smooth
affine algebraic variety over $\K$.
There is a canonical bijection
\[ \opn{quant} : \frac{ \{ \tup{formal Poisson brackets on 
$\mcal{O}_{X}[[\hbar]]$} \} }{\tup{ gauge equivalence}} 
\iso
\frac{ \{ \tup{star products on $\mcal{O}_X[[\hbar]]$} \} }
{\tup{ gauge equivalence}} .
 \] 
which is a quantization as defined above. 
\end{thm} 

By ``canonical'' we mean that this quantization map 
commutes with \'etale morphisms $X' \to X$ (and in particular with
automorphisms of $X$).

Actually our result in \cite{Ye1} is stronger -- it holds for a wider
class of varieties, not just affine varieties. 
However all these cases are subsumed in Corollary \ref{cor_new3} below. 

On the other hand, it might be good to remark that when writing \cite{Ye1} we
did not know Theorem \ref{thm:diff}, and hence we only considered differential
star products in that paper. Now we know that there is no difference (up to
gauge equivalence), so Theorem \ref{thm1} is correct as stated. In the context
of complex manifolds  this issue is still open: it is not known if every star
product is gauge equivalent to a differential one! See \cite{KS2}.

It is important to note that even the affine case of Theorem \ref{thm1} is a
``global result''. In this context ``local'' means a sufficiently small affine
open set $U \subset X$ that admits an \'etale coordinate system, namely an 
\'etale morphism $U \to \mbf{A}^n_{\K}$. This is totally analogous to the case
of $\mrm{C}^{\infty}$ manifolds studied by Kontsevich, where ``local'' meant an
open set in the manifold $\mbf{A}^n_{\R} = \R^n$.

Theorem \ref{thm1} is a consequence of the following more general
result.

\begin{thm}[\cite{Ye1}] \label{thm2} 
Let $\K$ be a field containing $\mbb{R}$, and let $X$ be a
smooth algebraic variety over $\K$. There is a diagram 
\[ \UseTips \xymatrix @C=5ex @R=3.5ex {
\mcal{T}^{}_{\mrm{poly}, X} 
\ar[d]
& \mcal{D}^{}_{\mrm{poly}, X}
\ar[d] \\
\opn{Mix}_{\bsym{U}}(\mcal{T}^{}_{\mrm{poly}, X} ) 
\ar[r]^{\Psi_{\bsym{\sigma}}}
&
\opn{Mix}_{\bsym{U}}(\mcal{D}^{}_{\mrm{poly}, X} )
} \] 
where:
\begin{itemize}
\item $\opn{Mix}_{\bsym{U}}(\mcal{T}^{}_{\mrm{poly}, X})$
and
$\opn{Mix}_{\bsym{U}}(\mcal{D}^{}_{\mrm{poly}, X})$
are sheaves of DG Lie algebras on $X$, called {\em mixed
resolutions};
\item the vertical arrows are DG Lie algebra quasi-isomorphisms;
\item and the horizontal arrow $\Psi_{\bsym{\sigma}}$
is an $\mrm{L}_{\infty}$ quasi-isomorphism. 
\end{itemize}
\end{thm}

The mixed resolutions combine the {\em commutative \v{C}ech 
resolution} associated to a sufficiently refined affine open covering $\bsym{U}$
of $X$, and the {\em Grothendieck sheaf of jets}. 
An {\em $\mrm{L}_{\infty}$ quasi-isomorphism} is a
generalization of a DG Lie algebra quasi-iso\-mor\-phism. 
The $\mrm{L}_{\infty}$ quasi-isomorphism $\Psi_{\bsym{\sigma}}$ depends on some
choices; but the dependence on these choices and on the covering 
$\bsym{U}$ disappears in homotopy. 
Theorem \ref{thm2} is proved using the {\em Formality Theorem}
of Kontsevich \cite{Ko1} and {\em formal geometry}.
More on the proof of Theorem \ref{thm2} in Appendices B and C.
A somewhat different approach to Theorem \ref{thm2} can be found in Van den
Bergh's paper \cite{VdB}.

\section{Twisted Deformations of Algebraic Varieties}
\label{sec:twdefs}

What can be done in general, when the variety $X$ is not affine? Can
we still make use of Theorem \ref{thm2}? 

In the paper \cite{Ko3} Kontsevich suggests that in general the 
deformation quantization of a Poisson bracket might have to be a 
{\em stack of algebroids}. This is a generalization of
the notion of sheaf of algebras. 

Actually stacks of algebroids appeared earlier, 
under the name {\em sheaves of twisted modules},
in the work of Kashiwara \cite{Ka}. See also \cite{DP}, \cite{PS},
\cite{KS1}, \cite{KS2}. 

I will use the term {\em twisted associative deformation},
and present an approach that treats the Poisson case as well. 
A similar point of view is taken in \cite{BGNT}. 

Here I will explain only a naive definition of twisted 
deformations. A more sophisticated definition, involving gerbes, 
may be found in Appendix D.
The fact that the two definitions agree follows from our work on
central extensions of gerbes and obstructions classes \cite{Ye5}.

Let $U \subset X$ be an affine open set, and let
$C := \Gamma(U, \mcal{O}_X)$. 
Suppose $A$ is an associative or Poisson deformation of the
$\K$-algebra $C$.
One may assume that 
$A = C[[\hbar]]$, and it is either endowed with a Poisson bracket
$\{ -,- \}$, or with a star product $\star$.  

In either case $A$ becomes a pronilpotent Lie algebra, and 
$\hbar A$ is a Lie subalgebra. In the Poisson
case the Lie bracket is $\{ -,- \}$, and in the associative case the
Lie bracket is the commutator
\[ [a,b] := a \star b - b \star a . \]
Let us denote the corresponding pronilpotent group by
\[ \opn{IG}(A) := \opn{exp}(\hbar A) , \]
and call it the {\em group of inner gauge transformations}
of $A$.

The group $\opn{IG}(A)$ acts on the deformation $A$ by gauge
equivalences. We denote this action by $\opn{ig}$. 
In the Poisson case the gauge
transformation $\opn{ig}(g)$, for $g \in \opn{IG}(A)$,
can be viewed as a formal hamiltonian flow. 
In the associative case the intrinsic exponential function
\[ \exp(a) =  \sum_{i \geq 0} \smfrac{1}{i!} \,
\underset{i}{\underbrace{ a \star \cdots \star a }} , \] 
for $a \in \hbar A$, allows us to identify the group $\opn{IG}(A)$
with the multiplicative subgroup 
\[ \{ g \in A \mid g \equiv 1 \opn{mod} \hbar \} . \]
Under this identification the operation $\opn{ig}(g)$ is just
conjugation by the invertible element $g$. 

The above can be sheafified: to a deformation $\mcal{A}$ of
$\mcal{O}_X$ we
assign the sheaf of groups $\opn{IG}(\mcal{A})$, etc.

Let us fix an affine open covering
$\{ U_0, \ldots, U_m \}$ of $X$. 
We write
\[ U_{i,j, \ldots} := U_i \cap U_j \cap \cdots . \]

\begin{dfn} \label{dfn1} 
A {\em twisted associative \tup{(}resp.\ Poisson\tup{)}
deformation} $\bsym{\mcal{A}}$ of $\mcal{O}_X$
consists of the following data: 
\begin{enumerate}
\item For any $i$, a deformation $\mcal{A}_i$ of $\mcal{O}_{U_i}$. 
\item For any $i < j$, a gauge equivalence
\[ g_{i,j} : \mcal{A}_i|_{U_{i,j}} \to \mcal{A}_j|_{U_{i,j}} .\] 
\item For any $i < j <k$, an element
\[ a_{i,j,k} \in \Gamma \bigl( U_{i,j,k}, \opn{IG}(\mcal{A}_i) \bigr)
. \]
\end{enumerate} 
The conditions are: 
\begin{enumerate}
\rmitem{i} For any $i <j <k$ one has
\[ g_{i,k}^{-1} \circ g_{j,k} \circ g_{i,j} = \opn{ig}(a_{i,j,k}^{-1})
. \]
\rmitem{ii} For any $i <j < k < l$ one has
\[ a_{i,j,l}^{-1} \cdot a_{i,k,l} \cdot a_{i,j,k} = 
g_{i,j}^{-1}(a_{j,k,l}^{}) . \]
\end{enumerate}
\end{dfn}

Condition (i) says that the $2$-cochain $\{ a_{i,j,k} \}$ 
measures the failure of the $1$-cochain $\{ g_{i,j} \}$ to be a
cocycle. This 
tells us whether the collection $\{ \mcal{A}_i \}$ of local
deformations can be glued into a global deformation of $\mcal{O}_X$. 

Condition (ii) -- usually called the tetrahedron equation -- says that
the $2$-cochain $\{ a_{i,j,k} \}$ satisfies a twisted cocycle condition. 

See Figure 1 for an illustration. 

\begin{figure}
\includegraphics[scale=0.4]{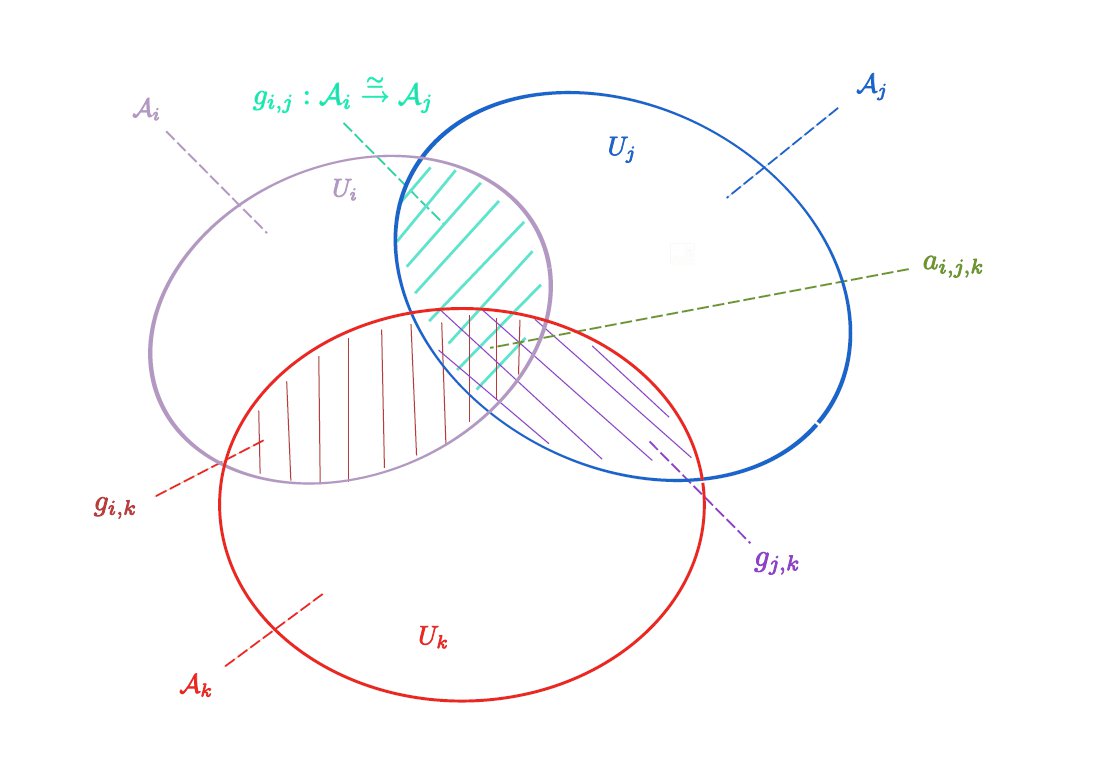}
\caption{}
\end{figure}

\begin{exa} \label{exa:20}
If $\mcal{A}$ is a usual deformation of $\mcal{O}_X$, then we obtain
a twisted deformation $\bsym{\mcal{A}}$ by taking
$\mcal{A}_i : = \mcal{A}|_{U_i}$, $g_{i,j} := \bsym{1}$ and
$a_{i,j,k} := 1$.

In this way we can view usual deformations as twisted deformations. 
\end{exa}

There is a notion of {\em twisted gauge equivalence} 
$\bsym{\mcal{A}} \to \bsym{\mcal{A}'}$ between twisted
associative (resp.\ Poisson) deformations of $\mcal{O}_X$. 

Let $\bsym{\mcal{A}}$ be a twisted deformation.
We say $\bsym{\mcal{A}}$ is {\em really twisted} if
it is not twisted-equivalent to any usual deformation $\mcal{A}'$. 

\begin{exa} 
It is easy to construct an example of a commutative twisted associative
deformation $\bsym{\mcal{A}}$ of $\mcal{O}_X$ that is really
twisted. Take an algebraic variety $X$ with nonzero cohomology
class $c \in \mrm{H}^2(X, \mcal{O}_X)$. 
Let $\bsym{U}$ be an affine open covering of $X$, and let
$\{ c_{i, j, k} \}$ be a \v{C}ech $2$-cocycle representing $c$
on this covering. Now consider the data
$\bigl( \{ \mcal{A}_i \}, \{ g_{i, j} \}, \{ a_{i, j, k} \} \bigr)$
with $\mcal{A}_i := \mcal{O}_{U_i}[[\hbar]]$, 
$g_{i, j} := \bsym{1}$ and
\[ a_{i, j, k} := \opn{exp}(\hbar c_{i, j, k}) . \]
This twisted deformation $\bsym{\mcal{A}}$ has obstruction
class $c$ in the first order central extension (in the sense of \cite{Ye6}).
This implies that $\bsym{\mcal{A}}$ is really twisted. 
\end{exa}

\begin{rem} \label{rem1}
For a twisted associative deformation $\bsym{\mcal{A}}$
there is a well defined abelian category 
$\cat{Coh} \bsym{\mcal{A}}$ of ``coherent left
$\bsym{\mcal{A}}$-modules'', which is a deformation of the 
abelian category $\cat{Coh} \mcal{O}_X$. See the work of Lowen and
Van den Bergh \cite{LV}.

Indeed, there is a
geometric Morita theory, which says that twisted associative
deformations of $\mcal{O}_X$ are the same as deformations of 
$\cat{Coh} \mcal{O}_X$. This is explained in the last chapter of the book
\cite{KS1} by Kashiwara and Schapira.

We do not know of a similar interpretation of twisted Poisson
deformations. 
\end{rem}

\section{Twisted Deformation Quantization}

Just as in the case of usual deformations, given a 
twisted (associative or Poisson) deformation $\bsym{\mcal{A}}$ 
of $\mcal{O}_X$, we can define the first order bracket 
$\{ -,- \}_{\bsym{\mcal{A}}}$ on $\mcal{O}_X$. The first order bracket is
invariant under twisted gauge equivalence.

Let $\bsym{\mcal{A}}$ be a twisted Poisson deformation, and 
let $\bsym{\mcal{B}}$ be a twisted associative deformation. We say 
that $\bsym{\mcal{B}}$ is a {\em twisted quantization} of
$\bsym{\mcal{A}}$ if 
\[ \{ -,- \}_{\bsym{\mcal{B}}} = \{ -,- \}_{\bsym{\mcal{A}}} . \]

The next theorem is influenced by ideas of Kontsevich (from 
\cite{Ko3} and private communications).

\begin{thm}[\cite{Ye5}]  \label{thm_new}
Let $\K$ be a field containing $\mbb{R}$, and let $X$ be a
smooth algebraic variety over $\K$.
There is a canonical bijection 
\[ \begin{aligned}
& \opn{tw{.}quant} : \frac{ \{ \tup{twisted Poisson deformations of 
$\mcal{O}_{X}$} \} }
{\tup{ twisted gauge equivalence}} \\
&   \hspace{3em} \iso \ 
\frac{ \{ \tup{twisted associative deformations of 
$\mcal{O}_X$} \} }{\tup{ twisted gauge equivalence}} ,
\end{aligned} \] 
called the {\em twisted quantization map}, which is a twisted quantization in
the sense above.
\end{thm} 

As before, by ``canonical'' we mean that the twisted quantization map
commutes with \'etale morphisms $X' \to X$.

The proof of Theorem \ref{thm_new} relies on a rather complicated
calculations involving Maurer-Cartan equations in cosimplicial DG Lie
algebras, and cosimplicial $2$-groupoids. See Appendices C and D. 

The results on pronilpotent gerbes from \cite{Ye6} tell us that 
when $\mrm{H}^1(X, \mcal{O}_{X}) = 0$ and 
$\mrm{H}^2(X, \mcal{O}_{X}) = 0$, there are no really twisted deformations. 
Thus Theorem \ref{thm_new} implies:

\begin{cor}  \label{cor_new}
Assume the cohomology groups 
$\mrm{H}^1(X, \mcal{O}_{X})$ and $\mrm{H}^2(X, \mcal{O}_{X})$ vanish.
Then the function $\opn{tw{.}quant}$ of the theorem gives a bijection
\[ \opn{quant} : \frac{ \{ \tup{Poisson deformations of 
$\mcal{O}_{X}$} \} }
{\tup{ gauge equivalence}} 
\iso
\frac{ \{ \tup{associative deformations of 
$\mcal{O}_X$} \} }{\tup{ gauge equivalence}} .
\] 
\end{cor} 

As we already mentioned, if $\mrm{H}^1(X, \mcal{T}_{X}) = 0$, then all
twisted Poisson deformations are sheaf theoretically trivial. 
Likewise, if $\mrm{H}^1(X, \mcal{D}_{X}) = 0$, then all
twisted associative deformations are sheaf theoretically trivial. 
Hence:

\begin{cor}  \label{cor_new3}
Assume the cohomology groups $\mrm{H}^1(X, \mcal{O}_{X})$,
$\mrm{H}^2(X, \mcal{O}_{X})$, \linebreak
$\mrm{H}^1(X, \mcal{T}_{X})$ and
$\mrm{H}^1(X, \mcal{D}_{X})$
vanish. Then the function $\opn{tw{.}quant}$ of the theorem gives a bijection
\[ \opn{quant} : \frac{ \{ \tup{formal Poisson brackets on 
$\mcal{O}_{X}[[\hbar]]$} \} }{\tup{ gauge equivalence}} 
\iso
\frac{ \{ \tup{star products on $\mcal{O}_X[[\hbar]]$} \} }
{\tup{ gauge equivalence}} .
 \] 
%
\end{cor} 

The conditions of Corollary \ref{cor_new3} are satisfied, for instance, when
$X = \mbf{P}^n_{\K}$ ($n$-dimensional projective space over $\K$). 

Let me finish this survey with a question. 
Given a variety $X$, with Poisson bracket $\{-,-\}_1$ on 
$\mcal{O}_X$, we can form the Poisson deformation 
$\mcal{A} := \mcal{O}_X[[\hbar]]$ with bracket $\hbar \{-,-\}_1$,
as in Example \ref{exa10}. 
By viewing $\mcal{A}$ as a twisted Poisson deformation
(cf.\ Example \ref{exa:20}), and applying
Theorem \ref{thm_new}, we get a twisted associative deformation 
$\bsym{\mcal{B}} := \opn{quant}(\mcal{A})$.

\begin{que}
Does there exist a variety $X$, with a symplectic
Poisson bracket $\{-,-\}_1$, such that the corresponding twisted
associative deformation $\bsym{\mcal{B}}$ is really twisted?
\end{que}

My feeling is that the answer is positive. 
And moreover, an example should be when $X$ is any Calabi-Yau surface,
and $\{-,-\}_1$ is any nonzero Poisson bracket on $X$.
It might be possible to settle this question with an explicit (yet very hard)
calculation, since there are explicit formulas for the twisted quantization
map.

\appendix

\section{DG Lie Algebras and Deformations}

The first time DG (differential graded) Lie algebras were used to study
problems in deformation theory was in the paper \cite{SS} of Schlessinger and
Stasheff (1977). In 1986 Deligne (in a letter to Millson, see \cite{GM})
formulated the idea that ``in characteristic zero, a deformation problem is
controlled by a differential graded Lie algebra, with quasi-isomorphic
differential graded Lie algebras giving the same deformation theory''.
This idea will be explained in this appendix, for associative and
Poisson deformations of smooth affine algebraic varieties. More details can be
found in the papers \cite{Ye1, Ye5}. 

Recall that a DG Lie algebra is a graded $\K$-module 
$\mfrak{g} = \boplus_{i \in \mbb{Z}} \mfrak{g}^{i}$,
with a bracket $[-, -]$ of degree $0$
satisfying the graded version of the Lie algebra identities, 
together with a graded derivation $\mrm{d}$ 
of degree $1$ and square $0$.

\begin{dfn}
A DG Lie algebra $\mfrak{g} = \boplus_{i \in \mbb{Z}} \mfrak{g}^{i}$
is said to be of {\em quantum type} if $\g^i = 0$ for all $i < -1$.
\end{dfn}

Given a DG Lie algebra $\mfrak{g}$, let us define a new DG Lie algebra 
\[ \mfrak{g}[[\hbar]]^+ := 
\bigoplus_{p \in \mbb{Z}} \, \hbar \mfrak{g}^{p}[[\hbar]] \subset
\bigoplus_{p \in \mbb{Z}} \, \mfrak{g}^{p}[[\hbar]] , \] 
in which $\hbar$ is central.

The {\em Maurer-Cartan equation} in 
$\mfrak{g}[[\hbar]]^+$ is
\[ \mrm{d}(\omega) + \smfrac{1}{2} [\omega, \omega] = 0 \]
for an element
\[ \omega = \sum_{j = 1}^{\infty} \omega_j \hbar^j
\in \mfrak{g}^{1}[[\hbar]]^+ . \]
The set of solutions of this equation is denoted by
$\opn{MC}(\mfrak{g}[[\hbar]]^+)$.

Let $\opn{exp}(\mfrak{g}^{0}[[\hbar]]^+)$
be the pronilpotent group associated to the pronilpotent Lie
algebra $\mfrak{g}^{0}[[\hbar]]^+$. 
There is an action of the group 
$\opn{exp}(\mfrak{g}^{0}[[\hbar]]^+)$ on 
$\mfrak{g}^{1}[[\hbar]]^+$, and this action preserves the subset 
$\opn{MC}(\mfrak{g}[[\hbar]]^+)$. One defines
\begin{equation} \label{eqn:MC}
\ol{\opn{MC}}(\mfrak{g}[[\hbar]]^+) := 
\frac{ \opn{MC}(\mfrak{g}[[\hbar]]^+) }{ 
\opn{exp}(\mfrak{g}^{0}[[\hbar]]^+) } \ ,
\end{equation}
the quotient set for this action.

Let us return to our deformation problem, where $X$ is a smooth 
algebraic variety over $\K$. Take an affine open set 
$U \subset X$, and let $C := \Gamma(U, \mcal{O}_X)$. 

One can show that any Poisson (resp.\ associative) deformation of $C$
is isomorphic to $C[[\hbar]]$ as $\K[[\hbar]]$-algebra (resp.\
$\K[[\hbar]]$-module) augmented to $C$. Thus it suffices to understand formal 
Poisson brackets (resp.\ star products) on $C[[\hbar]]$.

Let $\mcal{T}_C$ denote the module of derivations of $C$.  
For $p \geq -1$ define
\[ \mcal{T}^{p}_{\mrm{poly}}(C) := 
\bwedge_C^{p + 1} \mcal{T}_C . \]
So $\mcal{T}^{-1}_{\mrm{poly}}(C) = C$, 
$\mcal{T}^{0}_{\mrm{poly}}(C) = \mcal{T}_C$ and
$\mcal{T}^{1}_{\mrm{poly}}(C) = \bwedge_C^{2} \mcal{T}_C$.
The direct sum
\[ \mcal{T}_{\mrm{poly}}(C) := \bigoplus_{p \geq -1}
\mcal{T}^{p}_{\mrm{poly}}(C) \]
is a DG Lie algebra, called the algebra of {\em poly derivations} 
of $C$. The Lie bracket is the Schouten-Nijenhuis bracket, and the 
differential is $0$.

An element
\[ \omega = \sum_{j = 1}^{\infty} \omega_j \hbar^j
 \in \mcal{T}^{1}_{\mrm{poly}}(C)[[\hbar]]^+ \]
determines a skew-symmetric $\K[[\hbar]]$-bilinear biderivation of $C[[\hbar]]$,
namely  
\begin{equation}
\{ f, g \} := \sum_{j = 1}^{\infty} \omega_j(f, g) \hbar^j 
\end{equation}
for $f, g \in C$. A calculation shows
that $\omega$ satisfies the Maurer-Cartan equation if and only if
the corresponding biderivation $\{ -,- \}$ satisfies the Jacobi identity, i.e.\
it is a formal Poisson bracket on $C[[\hbar]]$.
Moreover, the action of group 
$\opn{exp} \bigl( \mcal{T}^0_{\mrm{poly}}(C)[[\hbar]]^+ \bigr)$
corresponds to gauge equivalences between formal Poisson brackets.
In this sense {\em $\mcal{T}_{\mrm{poly}}(C)$ controls Poisson 
deformations of $C$}. 

The second DG Lie algebra in this picture is that of the 
{\em poly differential operators}. 
A function $\phi : C^{p+1} \to C$ is called a poly differential
operator of order $\leq i$ if it is a differential operator of order $\leq i$ in
each of its arguments. The set of all such operators is denoted by
$\mcal{D}^{p}_{\mrm{poly}}(C)$.
Thus $\mcal{D}^{-1}_{\mrm{poly}}(C) = C$ and
$\mcal{D}^{0}_{\mrm{poly}}(C) = \mcal{D}(C)$, the ring of differential
operators. The direct sum 
\[ \mcal{D}_{\mrm{poly}}(C) := \bigoplus_{p \geq -1}
\mcal{D}^{p}_{\mrm{poly}}(C) \]
is a sub DG Lie algebra of the shifted Hochschild cochain complex
(with its Gerstenhaber Lie bracket). 

A solution
$\omega = \sum_{j = 1}^{\infty} \omega_j \hbar^j$
of the Maurer-Cartan equation in 
$\mcal{D}_{\mrm{poly}}(C)[[\hbar]]^+$ 
determines a star product on $C[[\hbar]]$, by the formula
\begin{equation} \label{eqn:beta}
f \star g := f g + 
\sum_{j = 1}^{\infty} \omega_j(f, g) \hbar^j
\end{equation}
for $f, g \in C$. 
Such a star product is called {\em differential} (as opposed to an ordinary star
product, for which the coefficients $\omega_j$ are just $\K$-bilinear; cf.\
equation (\ref{eqn:omega2})). The group 
$\opn{exp} \bigl( \mcal{D}^{0}_{\mrm{poly}}(C)[[\hbar]]^+ \bigr)$
is the group of differential gauge equivalences between differential star
products. Thus {\em the DG Lie algebra $\mcal{D}_{\mrm{poly}}(C)$ 
controls differential star products}.

Fortunately we have this result:

\begin{thm}[\cite{Ye5}] \label{thm:diff}
Let $C$ be a smooth $\K$-algebra. 
\begin{enumerate}
\item Any star product $\star$ on $C[[\hbar]]$ is gauge equivalent to a
differential star product $\star'$.

\item Let $\star$ and $\star'$ be differential star products on $C[[\hbar]]$, 
and let $g : \star \to \star'$ be a gauge equivalence. Then $g$ is a
differential gauge equivalence.
\end{enumerate}
\end{thm}

The proof of Theorem \ref{thm:diff} goes like this: it is well known that star
products on $C[[\hbar]]$ are controlled by the shifted Hochschild cochain
complex; this was essentially shown already in \cite{Ge}. We know that the
inclusion of the DG Lie algebra $\mcal{D}_{\mrm{poly}}(C)$ into the 
shifted Hochschild cochain complex is quasi-isomorphism (cf.\ \cite{Ye3}). 
This, together with the Equivalence Theorem (see below), implies
part (1) of Theorem \ref{thm:diff}. Part (2) is a direct calculation. 

The upshot is that the {\em DG Lie algebra 
$\mcal{D}_{\mrm{poly}}(C)$ in fact controls all star products} on $C[[\hbar]]$.

Both $\mcal{T}_{\mrm{poly}}(C)$ and $\mcal{D}_{\mrm{poly}}(C)$ are quantum type
DG Lie algebras (this is the reason for the name!).

Geometrically, there are sheaves of DG Lie algebras
$\mcal{T}_{\mrm{poly}, X}$ and $\mcal{D}_{\mrm{poly}, X}$
on $X$, that are quasi-coherent as $\mcal{O}_X$-modules. 
For any affine open set $U$ as above we have
\[ \Gamma(U, \mcal{T}_{\mrm{poly}, X}) = \mcal{T}_{\mrm{poly}}(C) , 
\]
and likewise for $\mcal{D}_{\mrm{poly}}$.

In order to control global deformations one has to resort to some 
kind of resolution of these sheaves of DG Lie algebras, such as the 
mixed resolutions mentioned in Theorem \ref{thm2}.

\begin{rem}
The Hochschild cochain complex is not
functorial in $C$ at all. On the other hand the complex 
$\mcal{D}^{}_{\mrm{poly}}(C)$ is functorial with respect to \'etale
homomorphisms $C \to C'$. 
\end{rem}

\begin{rem}
The product $\star$ in equation (\ref{eqn:beta}) could fail to be unital. 
To ensure that $\star$ has $1 \in C$ as unit, we have to take 
$\omega \in \opn{MC}(\mcal{D}^{\mrm{nor}}_{\mrm{poly}}(C)[[\hbar]]^+)$.
Here $\mcal{D}^{\mrm{nor}}_{\mrm{poly}}(C)$ is the algebra of {\em normalized
poly differential operators}. However, since the inclusion 
$\mcal{D}^{\mrm{nor}}_{\mrm{poly}}(C) \to \mcal{D}^{}_{\mrm{poly}}(C)$ 
is a quasi-isomorphism, the sets $\ol{\opn{MC}}(-)$ are the same. 
For this reason we can allow ourselves to neglect the distinction between
$\mcal{D}^{\mrm{nor}}_{\mrm{poly}}(C)$ and 
$\mcal{D}_{\mrm{poly}}(C)$ in this survey.
\end{rem}

\section{The Universal Quantization Map}

Let $C$ be a smooth $\K$-algebra. 
There is a canonical homomorphism of complexes of $\K$-modules
\[ \Upsilon_1 : \mcal{T}_{\mrm{poly}}(C) \to
\mcal{D}_{\mrm{poly}}(C) \]
given by
\[ \Upsilon_1(\partial_1 \wedge \cdots \wedge \partial_k)
(f_1, \ldots, f_k) := 
\smfrac{1}{k!} \sum_{\sigma \in S_k} \opn{sgn}(\sigma)
\partial_{\sigma(1)}(f_1) \cdots \partial_{\sigma(k)}(f_k) \]
for $f_i \in C$ and $\partial_i \in \mcal{T}_C$. 
Here $S_k$ is the group of permutations.
The homomorphism $\Upsilon_1$ is called the {\em antisymmetrization map}
or the {\em HKR map}, the latter because of its similarity to the famous
Hochschild-Kostant-Rosenberg Theorem. 
It is known that $\Upsilon_1$ is a quasi-isomorphism --
see \cite{Ko1} for the $\mrm{C}^{\infty}$ case,
and \cite{Ye1} for the algebraic case  --
and it induces an isomorphism of graded Lie algebras in cohomology. 
But $\Upsilon_1$ is not a DG Lie algebra homomorphism!

Here is the major result in the area of deformation quantization:

\begin{thm}[Kontsevich Formality Theorem, \cite{Ko1}, 1997] \label{thm:forma}
Let  $C :=$ \linebreak $\K[t_1, \ldots, t_n]$, the polynomial ring over $\K$,
and assume $\mbb{R} \subset \K$. Then $\Upsilon_1$ extends to an 
$\mrm{L}_{\infty}$ quasi-isomorphism
\[ \Upsilon = \{ \Upsilon_j \}_{j = 1}^{\infty} : 
\mcal{T}_{\mrm{poly}}(C) \to
\mcal{D}_{\mrm{poly}}(C) . \]
In other words, $\Upsilon_1$ is a DG Lie algebra 
quasi-isomorphism, up to specified higher homotopies 
$\Upsilon_2, \Upsilon_3, \ldots$. 
Each of the functions $\Upsilon_j$ is a poly differential operator, and is
invariant under linear change of coordinates in $C$.
\end{thm}

There is an induced $\K[[\hbar]]$-multilinear $\mrm{L}_{\infty}$
quasi-isomorphism
\begin{equation} \label{eqn:22}
\Upsilon : \mcal{T}_{\mrm{poly}}(C)[[\hbar]]^+ \to
\mcal{D}_{\mrm{poly}}(C)[[\hbar]]^+ . 
\end{equation}

The next result was known in the nilpotent case, namely for artinian parameter
algebras as in Remark \ref{rem:parameter}, at least since \cite {Ko1}.
The complete case was only proved recently \cite[Theorem 0.4]{Ye7}. 

\begin{thm}[Equivalence Theorem] \label{thm:equiv}
Let $\Phi : \g \to \h$ be an $\mrm{L}_{\infty}$ quasi-iso\-morphism between DG
Lie algebras. Then there is an induced bijection 
\[ \ol{\opn{MC}}(\Phi) : \ol{\opn{MC}}(\g[[\hbar]]^+) \to 
\ol{\opn{MC}}(\h[[\hbar]]^+) , \]
functorial in $\Phi$, with an explicit formula.
\end{thm}

Combining Theorems \ref{thm:forma} and \ref{thm:equiv} we find that there
is a canonical bijection
\[ \ol{\opn{MC}}(\Upsilon) : 
\ol{\opn{MC}} \bigl( \mcal{T}_{\mrm{poly}}(C)[[\hbar]]^+ \bigr) \iso
\ol{\opn{MC}} \bigl( \mcal{D}_{\mrm{poly}}(C)[[\hbar]]^+ \bigr) . \]
Therefore:

\begin{cor}
Assume  $\mbb{R} \subset \K$ and
$C = \K[t_1, \ldots, t_n]$. Then there is a canonical 
bijection of sets
\[ \opn{quant} : 
\frac{ \bigl\{ \tup{formal Poisson brackets on $C[[\hbar]]$}
\bigr\} }
{\tup{gauge equivalence}} \iso
\frac{ \bigl\{ \tup{star products on $C[[\hbar]]$} \bigr\} }
{\tup{gauge equivalence}} \]
preserving first order brackets. 
\end{cor}

\begin{rem}
The reason we require that $\K$ contains $\R$ is because the explicit formula
of Kontsevich for the higher homotopies $\Upsilon_j$ involves transcendental
real numbers. See discussion in \cite{Ko2}. More recent work (based on ideas of
Tamarkin, cf.\ \cite{CV2}) shows that this requirement can be avoided.
\end{rem}

\section{The $\mrm{L}_{\infty}$ quasi-isomorphism of the Level of 
Sheaves}

Here is an outline of the proof of Theorem \ref{thm2}.
Details can be found in the Erratum to \cite{Ye1}. 
We assume $\mbb{R} \subset \K$, and $X$ is a smooth $n$-dimensional
algebraic variety over $\K$.

A formal coordinate system at a closed point 
$x \in X$ is an isomorphism of $\K$-algebras
\[ \bsym{k}(x)[[\bsym{t}]] = 
\bsym{k}(x)[[t_1, \ldots, t_n]] \iso \what{\mcal{O}}_{X, x}
, \]
where $\bsym{k}(x)$ is the residue field, and 
$\what{\mcal{O}}_{X, x}$ is the complete local ring.
 
There is an infinite dimensional scheme 
$\opn{Coor} X$, called the {\em coordinate bundle}, with a projection 
$\pi : \opn{Coor} X \to X$, which is a moduli space for
formal coordinate systems. (In \cite{Ko1} the notation for 
$\opn{Coor} X$ is $X^{\mrm{coor}}$.)
In particular, for every closed 
point $x \in X$, the $\bsym{k}(x)$-rational points 
in the fiber $\pi^{-1}(x)$ stand in bijection to 
the set of formal coordinate systems at $x$.

To get an idea of how the scheme $\opn{Coor} X$ looks, let us note 
that 
$\opn{Coor} X = \opn{lim}_{\leftarrow} \opn{Coor}^i X$,
where each $\opn{Coor}^i X$ is the variety 
parameterizing formal coordinate systems up to order $i$. 

Let $\mcal{P}_X$ be the sheaf of principal
parts on $X$ (a.k.a.\ the jet sheaf). Recall that $\mcal{P}_X$
is the formal completion of 
$\mcal{O}_{X \times X}$ along the diagonal. As such, $\mcal{P}_X$ is a sheaf 
of rings on $X$, with two ring homomorphisms 
$\mcal{O}_X \to \mcal{P}_X$ (corresponding to the two projections 
$X \times X \to X$), which make $\mcal{P}_X$ into an 
$\mcal{O}_X$-bimodule. When we view $\mcal{P}_X$ as a left $\mcal{O}_X$-module,
it has the  {\em Grothendieck connection}
\[ \nabla : \mcal{P}_X \to \Omega^1_X \otimes_{\mcal{O}_X} \mcal{P}_X , \]
which is the completion of pullback, under the first projection 
$X \times X \to X$, of the standard connection 
$\d : \mcal{O}_X \to \Omega^1_X$ . The Grothendieck connection is flat, and its
kernel is
$\mcal{O}_X$ (coming from the second projection). 

Let us denote by $\pi^* \mcal{P}_X$ the pullback of the left 
$\mcal{O}_X$-module $\mcal{P}_X$; and let 
$\pi^{\what{*}} \mcal{P}_X$ be its adic completion. 
The universal property of 
$\pi : \opn{Coor} X \to X$ implies that there is a canonical isomorphism 
\begin{equation} \label{eqn:univ}
\mcal{O}_{\opn{Coor} X} \hatotimes{\K} \K[[\bsym{t}]]
\cong \pi^{\what{*}} \mcal{P}_X ,
\end{equation}
of sheaves of $\mcal{O}_{\opn{Coor} X}$-algebras on $\opn{Coor} X$. 

Like (\ref{eqn:univ}), there are canonical isomorphisms
\[ \mcal{O}_{\opn{Coor} X} \hatotimes{\K} 
\mcal{T}_{\mrm{poly}}(\K[[\bsym{t}]]) \cong
\pi^{\what{*}} (\mcal{P}_X \otimes_{\mcal{O}_X} 
\mcal{T}_{\mrm{poly}, X}) \]
and
\[ \mcal{O}_{\opn{Coor} X} \hatotimes{\K} 
\mcal{D}_{\mrm{poly}}(\K[[\bsym{t}]]) \cong 
\pi^{\what{*}} (\mcal{P}_X \otimes_{\mcal{O}_X} 
\mcal{D}_{\mrm{poly}, X}) \]
of graded Lie algebras on $\opn{Coor} X$. Let us write 
$\mcal{A} := \Omega_{\opn{Coor} X}$, the sheaf of differential forms on
$\opn{Coor} X$; so in particular 
$\mcal{A}^0 = \mcal{O}_{\opn{Coor} X}$.
Then, by applying $\mcal{A} \hatotimes{\mcal{A}^0} -$, we get canonical
graded Lie algebra isomorphisms 
\begin{equation} \label{eqn:50}
\mcal{A}  \hatotimes{\K} 
\mcal{T}_{\mrm{poly}}(\K[[\bsym{t}]]) \cong
\mcal{A} \hatotimes{\mcal{A}^0}
\pi^{\what{*}} (\mcal{P}_X \otimes_{\mcal{O}_X} 
\mcal{T}_{\mrm{poly}, X})
\end{equation}
and
\begin{equation} \label{eqn:51}
\mcal{A}  \hatotimes{\K} 
\mcal{D}_{\mrm{poly}}(\K[[\bsym{t}]]) \cong
\mcal{A} \hatotimes{\mcal{A}^0}
\pi^{\what{*}} (\mcal{P}_X \otimes_{\mcal{O}_X} 
\mcal{D}_{\mrm{poly}, X}) 
\end{equation}
on $\opn{Coor} X$. However the differentials do not match: to compensate for the
Grothen\-dieck connection on the right side, we have to add the differential 
$\opn{ad}(\omega)$ on the left side, where 
$\omega \in \mcal{A}^1 \hatotimes{\K} \mcal{T}^0_{\mrm{poly}}(\K[[\bsym{t}]])$
is a universal MC element.

Due to the Kontsevich Formality Theorem (Theorem \ref{thm:forma}) we obtain an
induced $\mcal{A}$-multilinear $\mrm{L}_{\infty}$ quasi-isomorphism
\[ \Upsilon : \mcal{A} \hatotimes{\K} 
\mcal{T}_{\mrm{poly}}(\K[[\bsym{t}]]) \to
\mcal{A} \hatotimes{\K} 
\mcal{D}_{\mrm{poly}}(\K[[\bsym{t}]])  \]
between sheaves of DG Lie algebras on $\opn{Coor} X$.
(This should be compared to the $\mrm{L}_{\infty}$ quasi-isomorphism 
(\ref{eqn:22}) in the local case.)
Using the isomorphisms (\ref{eqn:50}) and (\ref{eqn:51}), and twisting 
$\Upsilon$ by the element $\omega$ (in the sense of \cite[Theorem 3.2]{Ye3}), we
then
obtain an $\mcal{A}$-multilinear $\mrm{L}_{\infty}$ quasi-isomorphism
\[ \Upsilon_{\omega} : \mcal{A} \hatotimes{\mcal{A}^0}
\pi^{\what{*}} (\mcal{P}_X \otimes_{\mcal{O}_X} \mcal{T}_{\mrm{poly}, X}) \to 
\mcal{A} \hatotimes{\mcal{A}^0}
\pi^{\what{*}} (\mcal{P}_X \otimes_{\mcal{O}_X} \mcal{D}_{\mrm{poly}, X}) . \]

If we had a section $\sigma : X \to \opn{Coor} X$
then we could pull $\Upsilon_{\omega}$ down to an $\mrm{L}_{\infty}$ 
quasi-isomorphism on $X$. 
However usually there are no global sections of $\opn{Coor} X$, because of
topological obstructions.

The group $\opn{GL}_{n}$ acts on $\opn{Coor} X$ by linear change 
of coordinates. Let us define $\opn{LCC} X$ to be the quotient scheme
$\opn{Coor} X / \opn{GL}_{n}$. (``LCC'' stands for ``linear 
coordinate classes''.) So the projection 
$\pi : \opn{Coor} X \to X$ factors through $\opn{LCC} X$.
Recall that the universal 
quantization of Kontsevich is invariant under linear change 
of coordinates, namely under the action of the group 
$\opn{GL}_{n}$ and its Lie algebra. This implies that the $\mrm{L}_{\infty}$ 
morphism $\Upsilon_{\omega}$ descends to $\opn{LCC} X$; 
and hence it suffices to look for sections 
$\sigma : X \to \opn{LCC} X$ of the projection
$\opn{LCC} X \to X$.

In the $\mrm{C}^{\infty}$ context such global sections 
$\sigma : X \to \opn{LCC} X$ do exists (because the fibers of the 
bundle $\opn{LCC} X$ are contractible). But this is not the case 
in our algebraic situation; so we must use a trick.

Let $G$ be the group of $\K$-algebra automorphisms of 
$\K[[\bsym{t}]]$. So
$G \cong \opn{GL}_n \ltimes N$,
where $N$ is the subgroup of elements that act trivially modulo 
$(\bsym{t})^2$. The group $N$ is pro-unipotent. 
It turns out that $\opn{Coor} X$ is a $G$-torsor over $X$, locally trivial in
the Zariski topology. 

Suppose we are given a finite number of sections
\[ \sigma_0, \ldots, \sigma_q : U \to \opn{LCC} X \]
over some open set $U \subset X$. 
Using an averaging process for unipotent 
group actions \cite{Ye4}, we show that there exists a morphism 
\[ \bsym{\sigma} : \bsym{\Delta}^q_{\K} \times U  \to \opn{LCC} X \] 
which restricts to $\sigma_j$ on the $j$-th vertex of
$\bsym{\Delta}^q_{\K}$. Here $\bsym{\Delta}^q_{\K}$ is the 
$q$-dimensional geometric simplex.

\begin{figure}
\includegraphics[scale=0.4]{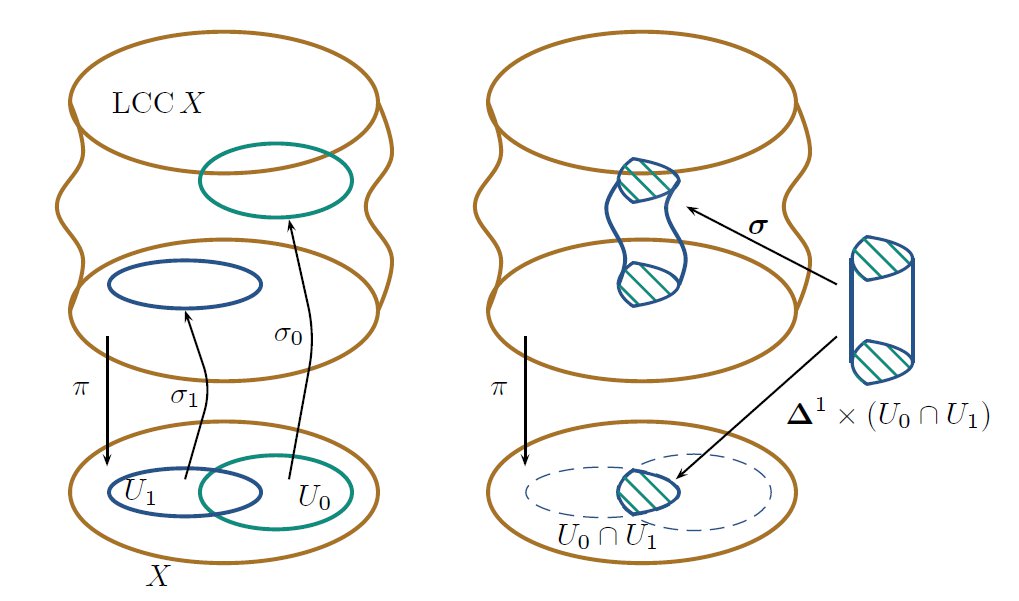}
\caption{Simplicial sections, $q = 1$. 
We start with sections over two open sets $U_0$ and $U_1$ 
in the left diagram; and we pass to a simplicial 
section $\bsym{\sigma}$ on the right.} 
\label{fig:1}
\end{figure}

Since sections exist locally, we can choose an open covering 
$\bsym{U} = \{ U_i \}$ of $X$, with sections 
$\sigma_i : U_i \to \opn{LCC} X$. For any
$i_0, \ldots, i_q$ we then obtain a morphism
\[ \bsym{\sigma} : \bsym{\Delta}^q_{\K} \times 
(U_{i_0} \cap \cdots \cap U_{i_q}) \to \opn{LCC} X . \]
(See Figure \ref{fig:1} for an illustration of the case $q = 1$.)
As $q$ varies we have a {\em simplicial section}
of the projection $\opn{LCC} X \to X$. Details are in \cite{Ye2}.

Another device we use is {\em mixed resolutions}. The mixed 
resolution $\opn{Mix}_{\bsym{U}}(\mcal{T}_{\mrm{poly}, X})$
is a sheaf of DG Lie algebras on $X$, depending on the covering
$\bsym{U}$, and equipped with a DG Lie algebra quasi-isomorphism 
\[ \mcal{T}_{\mrm{poly}, X} \to \opn{Mix}_{\bsym{U}}(\mcal{T}_{\mrm{poly}, X})
. \]
Likewise for $\mcal{D}_{\mrm{poly}, X}$. 

The simplicial section $\bsym{\sigma}$
allows us to pull the $\mrm{L}_{\infty}$ quasi-isomorphism
$\Upsilon_{\omega}$ down to $X$, giving  an $\mrm{L}_{\infty}$ quasi-isomorphism
\[ \Psi_{\bsym{\sigma}} :
\opn{Mix}_{\bsym{U}}(\mcal{T}_{\mrm{poly}, X}) \to
\opn{Mix}_{\bsym{U}}(\mcal{D}_{\mrm{poly}, X})  \]
between these sheaves of DG Lie algebras on $X$.
The dependence on $\bsym{U}$ and $\bsym{\sigma}$ disappears when we pass to
homotopy classes.

\section{Twisted Deformations and their Gauge Gerbes}

This section summarizes the main concepts of the paper \cite{Ye5}. 

Recall that a {\em groupoid} $\cat{G}$ is a category in which
all morphisms are invertible. We denote by $\cat{G}(i,j)$ the set of arrows from
the object $i$ to the object $j$. Note that $\cat{G}(i,i)$ is a group. %
Any element $g \in \cat{G}(i,j)$ defines a group isomorphism
\[ \opn{Ad}(g) : \cat{G}(i,i) \iso \cat{G}(j,j) , \]
namely 
$\opn{Ad}(g)(h) := g \circ h \circ g^{-1}$.

A {\em stack of groupoids} $\bsym{\mcal{G}}$ on $X$
is the geometrization of the notion of groupoid, 
in the same way that a sheaf of groups is the geometrization of the
notion of a group. 
Thus for any open set $U \subset X$ there is a groupoid
$\bsym{\mcal{G}}(U)$. 
And there are restriction functors
$\bsym{\mcal{G}}(U) \to \bsym{\mcal{G}}(V)$
for any inclusion $V \subset U$. 
These satisfy a rather complicated list of conditions. 
In particular, given any open set $U \subset X$ and any object 
$i \in \opn{Ob} \bsym{\mcal{G}}(U)$,
there is a sheaf of groups $\bsym{\mcal{G}}(i,i)$ on $U$. 
For details see \cite[Sections 1-2]{Ye6}, or
\cite{Gi, BM, KS1}. 

A stack of groupoids $\bsym{\mcal{G}}$ is called a {\em gerbe} if
it is locally nonempty and locally connected.

\begin{dfn} \label{dfn2}
Let $X$ be a smooth algebraic variety over $\K$. 
A {\em twisted associative \tup{(}resp.\ Poisson\tup{)}
deformation} $\bsym{\mcal{A}}$ of $\mcal{O}_X$ is 
the following data: 
\begin{enumerate}
\item A gerbe $\bsym{\mcal{G}}$ on $X$, called the
{\em gauge gerbe} of $\bsym{\mcal{A}}$. 
\item For any open set $U \subset X$ and 
$i \in \opn{Ob} \bsym{\mcal{G}}(U)$, an associative 
\tup{(}resp.\ Poisson\tup{)} deformation 
$\mcal{A}_i$ of $\mcal{O}_{U_i}$.
\end{enumerate} %
The conditions are: %
\begin{enumerate}
\rmitem{a} For any $i \in \opn{Ob} \bsym{\mcal{G}}(U)$, the sheaf of
groups
$\bsym{\mcal{G}}(i,i)$ coincides with 
$\opn{IG}(\mcal{A}_i)$, the sheaf of inner gauge transformations of
the
deformation $\mcal{A}_i$. %
\rmitem{b} For any $i \in \opn{Ob} \bsym{\mcal{G}}(U)$, any
$j \in \opn{Ob} \bsym{\mcal{G}}(V)$, any $W \subset U \cap V$ and any
$g \in \bsym{\mcal{G}}(W)(i,j)$, the isomorphism of sheaves of groups
\[ \opn{Ad}(g) : \mcal{G}(i,i)|_W \iso \mcal{G}(j,j)|_W \]
is induced from a (necessarily unique) gauge equivalence
\[ \mcal{A}_i|_W \iso \mcal{A}_j|_W . \]
\end{enumerate}
\end{dfn}

\begin{thm}[\cite{Ye6}]  \label{thm.last}
Definitions \tup{\ref{dfn1}} and \tup{\ref{dfn2}} are equivalent.
\end{thm} 

The proof relies on the fact that the gauge gerbe $\bsym{\mcal{G}}$ is
pronilpotent, and its abelian slices are coherent $\mcal{O}_X$-modules. Hence
for any affine open set $U$ the groupoid 
$\bsym{\mcal{G}}(U)$ is nonempty and connected.

Finally let us say a few words on the proof of Theorem \ref{thm_new}
in \cite{Ye5} \footnote{This is the proof in the new version of \cite{Ye5},
still in preparation as of July 2011.}.

Fix an affine open covering 
$\bsym{U} = \{ U_0, \ldots, U_m \}$ of $X$, such that for each $i$
there is an \'etale morphism 
$U_i \to \mbf{A}^n_{\K}$. 
Given a sheaf $\mcal{G}$ of DG Lie algebras on $X$, the \v{C}ech construction
gives a cosimplicial DG Lie algebra
$\opn{C}(\bsym{U}, \mcal{G}) = 
\{ \opn{C}^p(\bsym{U}, \mcal{G}) \}_{p \in \mbb{N}}$.
Thus from  Theorem \ref{thm2} we deduce that there is a diagram
\begin{equation} \label{eqn:last}
\UseTips \xymatrix @C=10ex @R=3.5ex {
\opn{C}(\bsym{U}, \mcal{T}^{}_{\mrm{poly}, X}) 
\ar[d]
& \opn{C}(\bsym{U}, \mcal{D}^{}_{\mrm{poly}, X})
\ar[d] \\
\opn{C}(\bsym{U}, \opn{Mix}_{\bsym{U}}(\mcal{T}^{}_{\mrm{poly}, X} ))
\ar[r]^{\opn{C}(\bsym{U}, \Psi_{\bsym{\sigma}})}
&
\opn{C}(\bsym{U}, \opn{Mix}_{\bsym{U}}(\mcal{D}^{}_{\mrm{poly}, X} ))
} 
\end{equation}
in which the objects are cosimplicial DG Lie algebras, the vertical arrows 
are cosimplicial DG Lie algebra homomorphisms, and the horizontal arrow is a
cosimplicial $\mrm{L}_{\infty}$ morphism. 
Moreover, in each cosimplicial dimension $p$ the diagram 
\[ \UseTips \xymatrix @C=10ex @R=3.5ex {
\opn{C}^p(\bsym{U}, \mcal{T}^{}_{\mrm{poly}, X}) 
\ar[d]
& \opn{C}^p(\bsym{U}, \mcal{D}^{}_{\mrm{poly}, X})
\ar[d] \\
\opn{C}^p(\bsym{U}, \opn{Mix}_{\bsym{U}}(\mcal{T}^{}_{\mrm{poly}, X} ))
\ar[r]^{\opn{C}^p(\bsym{U}, \Psi_{\bsym{\sigma}})}
&
\opn{C}^p(\bsym{U}, \opn{Mix}_{\bsym{U}}(\mcal{D}^{}_{\mrm{poly}, X} ))
} \] 
has objects that are quantum type DG Lie algebras, the vertical arrows 
are DG Lie algebra quasi-isomorphisms, and the horizontal arrow is an
$\mrm{L}_{\infty}$ quasi-isomorphism.

Now take any quantum type DG Lie algebra $\g$. To it we associate, in a
functorial way, the {\em Deligne $2$-groupoid} 
$\bsym{\cat{Del}}(\g[[\hbar]]^+)$, as in \cite{Ge}. 
If we are given a cosimplicial quantum type DG Lie algebra 
$\g = \{ \g^p \}_{p \in \mbb{N}}$, then there is  a cosimplicial
$2$-groupoid
\[ \bsym{\cat{Del}}(\g[[\hbar]]^+) := 
\bigl\{ \bsym{\cat{Del}}(\g^p[[\hbar]]^+) \bigr\}_{p \in \mbb{N}} . \]

Next let $\bsym{\cat{G}} = \{ \bsym{\cat{G}}^p \}_{p \in \mbb{N}}$ be any 
cosimplicial $2$-groupoid. To it we associate the set 
$\opn{CDD}(\bsym{\cat{G}})$ of {\em combinatorial descent data}, and the
quotient set by gauge equivalences $\ol{\opn{CDD}}(\bsym{\cat{G}})$.
The sets $\opn{CDD}(\bsym{\cat{G}})$ and $\ol{\opn{CDD}}(\bsym{\cat{G}})$
depend functorially on $\bsym{\cat{G}}$.

So for a cosimplicial quantum type DG Lie algebra 
$\g = \{ \g^p \}_{p \in \mbb{N}}$ we can associate the set
\[ \opn{LDD}(\g[[\hbar]]^+) := \opn{CDD} \bigl( \bsym{\cat{Del}}(\g[[\hbar]]^+)
\bigr) \]
of {\em Lie descent data}
\footnote{In the version of \cite{Ye5} dated August 2009, 
this is called ``additive descent data'', and denoted by 
$\opn{ADD}(\g[[\hbar]]^+)$.}, 
and its quotient set 
\[ \ol{\opn{LDD}}(\g[[\hbar]]^+) := \ol{\opn{CDD}} 
\bigl( \bsym{\cat{Del}}(\g[[\hbar]]^+) \bigr) . \]
See \cite[Section 3.2]{BGNT}, where $\opn{LDD}(\g[[\hbar]]^+)$
is called ``descent data for Deligne $2$-groupoids''.  

We have the following theorem, extending \cite[Proposition 3.3.1]{BGNT}.

\begin{thm}[\cite{Ye5}]
Let $\g$ and $\h$ be cosimplicial quantum type DG Lie algebras, and let 
$\Phi : \g \to \h$ be a cosimplicial $\mrm{L}_{\infty}$ quasi-isomorphism.
Then there is a bijection 
\[ \ol{\opn{LDD}}(\Phi) : \ol{\opn{LDD}}(\g[[\hbar]]^+) \iso
\ol{\opn{LDD}}(\h[[\hbar]]^+) , \]
depending functorially on $\Phi$, with an explicit formula.
\end{thm}

Applying this theorem to the diagram (\ref{eqn:last}) we obtain a canonical
diagram of bijections of sets
\[ \UseTips \xymatrix @C=8ex @R=3.5ex {
\ol{\opn{LDD}} \bigl( 
\opn{C}(\bsym{U}, \mcal{T}^{}_{\mrm{poly}, X})[[\hbar]]^+ \bigr)
\ar[d]
\ar@{-->}[r]^{\opn{tw{.}quant}}
& 
\ol{\opn{LDD}} \bigl( 
\opn{C}(\bsym{U}, \mcal{D}^{}_{\mrm{poly}, X})[[\hbar]]^+ \bigr)
\ar[d] 
\\
\ol{\opn{LDD}} \bigl( 
\opn{C}(\bsym{U},
\opn{Mix}_{\bsym{U}}(\mcal{T}^{}_{\mrm{poly}, X} ))
[[\hbar]]^+ \bigr)
\ar[r]
&
\ol{\opn{LDD}} \bigl( 
\opn{C}(\bsym{U}, \opn{Mix}_{\bsym{U}}(\mcal{D}^{}_{\mrm{poly}, X} ))
[[\hbar]]^+ \bigr) \ .
} \] 
Theorem \ref{thm.last} implies that twisted Poisson deformations of
$\mcal{O}_X$, modulo twisted gauge equivalence, correspond bijectively
to elements of  
\[ \ol{\opn{LDD}} \bigl( \opn{C}(\bsym{U}, \mcal{T}^{}_{\mrm{poly}, X})
[[\hbar]]^+ \bigr) . \]
Likewise, twisted associative deformations,
modulo twisted gauge equivalence, correspond bijectively to elements of
\[ \ol{\opn{LDD}} \bigl( \opn{C}(\bsym{U}, \mcal{D}^{}_{\mrm{poly}, X})
[[\hbar]]^+ \bigr) . \]
The resulting bijection
\[ \begin{aligned}
& \opn{tw{.}quant} : \frac{ \{ \tup{twisted Poisson deformations of 
$\mcal{O}_{X}$} \} }
{\tup{ twisted gauge equivalence}} \\
&  \hspace{3em} \iso \ 
\frac{ \{ \tup{twisted associative deformations of 
$\mcal{O}_X$} \} }{\tup{ twisted gauge equivalence}} 
\end{aligned} \] 
is canonical (independent of $\bsym{U}$ and $\bsym{\sigma}$).

\bibliographystyle{amsplain}

\begin{thebibliography}{BFFLS}
\bibitem[BFFLS]{BFFLS} F. Bayen, M. Flato, C. Fronsdal, 
    A. Lichnerowicz and D. Sternheimer,
    Deformation theory and quantization I, II, Ann.\ Physics
    {\bf 111} (1978), 61-110, 111-151.
\bibitem[BGNT]{BGNT} P. Bressler, A. Gorokhovsky, R. Nest and
	B. Tsygan, Deformation quantization of gerbes,
	Adv.\ Math.\ {\bf 214}, Issue 1 (2007), 230-266.
\bibitem[BK]{BK} R. Bezrukavnikov and D. Kaledin, 
    Fedosov quantization in algebraic context, 
    Moscow Math.\ J. {\bf 4}, no.\ 3 (2004), 559592.
\bibitem[BS]{BS} J. Baez and U. Schreiber, Higher Gauge Theory, 
  in: ``Categories in Algebra, Geometry and Mathematical Physics'',
  Contemporary Mathematics {\bf 431} (2007), AMS.
\bibitem[BM]{BM} L. Breen and W. Messing, 
    Differential geometry of gerbes,
    Adv.\ Math.\ {\bf 198} (2005), no.\ 2, 732-846.
\bibitem[CFT]{CFT} A. Cattaneo, G. Felder and L. Tomassini, 
    From local to global deformation quantization of Poisson 
    manifolds, Duke Math.\ J.\ \textbf{115} (2002), no.\ 2, 329-352.
\bibitem[CKTB]{CKTB} A. Cattaneo, B. Keller, C. Torossian and
    A. Bruguieres, 
    ``D\'eformation, Quantification, Th\'eory de Lie'',
    Panoramas et Synth\`eses {\bf 20} (2005), Soc.\ Math.\ France.
\bibitem[CV1]{CV1} D. Calaque and M. Van den Bergh,
	Hochschild cohomology and Atiyah classes,
  Advances in Mathematics {\bf 224} (2010), 1839-1889.
\bibitem[CV2]{CV2} D. Calaque and M. Van den Bergh,
    Global formality at the G-infinity level, 
  Mosc.\ Math.\ J. {\bf 10} (2010), 31-64, 271.
\bibitem[DL]{DL} M. De Wilde and P.B.A. Lecomte, 
    Existence of star-products and of formal deformations in 
    Poisson Lie algebra of arbitrary symplectic manifolds, Let.\ 
    Math.\ Phys.\ {\bf 7} (1983), 487-496.
\bibitem[DP]{DP}  A. D'Agnolo and P. Polesello,
    Stacks of twisted modules and integral transforms,
  in: ``Geometric aspects of Dwork theory'', Vol.\ I, II, Walter de Gruyter GmbH
  and Co.\ KG, Berlin (2004), 463-507.
\bibitem[Fe]{Fe} B. Fedosov, 
    A simple geometrical construction of deformation quantization,
    J. Differential Geom.\ {\bf 40} (1994), no.\ 2, 21-238. 
\bibitem[Gi]{Gi} J. Giraud, 
    ``Cohomologie non abelienne,'' 
    Grundlehren der Math.\ Wiss.\ {\bf 179}, Springer (1971).
\bibitem[Ge]{Ge} M. Gerstenhaber,
    On the deformation of rings and algebras,
    Ann.\ of Math.\ {\bf 79} (1964), 59-103.
\bibitem[GK]{GK} I.M. Gelfand and D.A. Kazhdan, 
    Some problems of differential geometry and the calculation of 
    cohomologies of Lie algebras of vector fields, 
    Soviet Math.\ Dokl.\ \textbf{12} (1971), no.\ 5, 1367-1370. 
\bibitem[GM]{GM} W.M. Goldman and J.J Millson, 
   The deformation theory of representations of fundamental groups of
      compact Kahler manifolds, Publ.\ Math.\ IHES
      {\bf 67} (1988), 43-96.
\bibitem[GS]{GS} V.W. Guillemin and S. Sternberg,
    ``Symplectic Techniques in Physics,'' 
    Cambridge University Press, Reprint 1990.
\bibitem[Hi]{Hi} V. Hinich, 
    Descent of Deligne groupoids,
    Internat.\ Math.\ Res.\ Notices {\bf 5} (1997), 223-239.
\bibitem[HY]{HY} R.\ H\"{u}bl and A.\ Yekutieli,
    Adelic Chern forms and applications,  {\em Amer.\ J. Math.}
    \textbf{121} (1999), 797-839.
\bibitem[Ka]{Ka} M. Kashiwara, 
    Quantization of contact manifolds, 
    Publ.\ Res.\ Inst.\ Math.\ Sci.\ {\bf 32} no.\ 1 (1996), 17.
\bibitem[Ko1]{Ko1} M. Kontsevich,
    Deformation quantization of Poisson manifolds,
    Lett.\ Math.\ Phys.\ {\bf 66} (2003), no.\ 3, 157-216.
\bibitem[Ko2]{Ko2} M. Kontsevich,
    Operads and Motives in deformation quantization, 
    Lett.\ Math.\ Phys.\ {\bf 48} (1999), 35-72.
\bibitem[Ko3]{Ko3} M. Kontsevich,
    Deformation quantization of algebraic varieties,
    Lett.\ Math.\ Phys.\ {\bf 56} (2001), no.\ 3, 271-294.
\bibitem[KS1]{KS1} M. Kashiwara and P. Schapira,
    ``Categories and Sheaves,'' Springer, 2006.
\bibitem[KS2]{KS2} M. Kashiwara and P. Schapira,
  Deformation quantization modules, arXiv:1003.3304.
\bibitem[Lo]{Lo} W. Lowen, 
  Algebroid prestacks and deformations of ringed spaces,
  Trans.\ Amer.\ Math.\ Soc.\ {\bf 360} (2008), 1631-1660. 
\bibitem[LV]{LV} W. Lowen and M. Van den Bergh,
    Deformation theory of abelian categories,
    Trans.\ AMS {\bf 358} (2006) no.\ 12, 5441 - 5483.
\bibitem[Mo]{Mo} I. Moerdijk, 
    Introduction to the Language of Stacks and Gerbes,
    eprint math.AG/0212266 at http://arxiv.org.
\bibitem[NT]{NT} R. Nest and B. Tsygan, 
    Deformations of symplectic Lie algebroids, 
    deformations of holomorphic symplectic structures, 
    and index theorems,  
    Asian J. Math.\ \textbf{5} (2001), no.\ 4, 599-635. 
\bibitem[PS]{PS} P. Polesello and P. Schapira, 
    Stacks of quantization-deformation modules on complex
    symplectic manifolds, 
    Intern.\ Math.\ Res.\ Notices. {\bf 2004} (2004), 2637-2664. 
\bibitem[SS]{SS} M. Schlessinger and J. Stasheff,
  Rational homotopy theory -- obstructions and deformations,
  In Proc.\ Conf.\ on Algebraic Topology, Vancouver, pages 7-31, 1977. 
  LMM 673.
\bibitem[VdB]{VdB} M. Van den Bergh,
    On global deformation quantization in the algebraic case,
    J. Algebra {\bf 315} (2007), 326-395.
\bibitem[Ye1]{Ye1} A.\ Yekutieli, 
    Deformation Quantization in Algebraic Geometry,
    Adv.\ Math.\ {\bf 198} (2005), 383-432.
    Erratum: Adv.\ Math.\ {\bf 217} (2008), 2897-2906.
\bibitem[Ye2]{Ye2} A.\ Yekutieli,
    Mixed Resolutions and Simplicial Sections,
    Israel J. Math.\ {\bf 162} (2007), 1-27.
\bibitem[Ye3]{Ye3} A.\ Yekutieli,
    Continuous and Twisted L-infinity Morphisms,
    J. Pure Appl.\ Algebra {\bf 207} (2006), 575-606.
\bibitem[Ye4]{Ye4} A.\ Yekutieli,
    An Averaging Process for Unipotent Group Actions,
    Representation Theory {\bf 10} (2006), 147-157.
\bibitem[Ye5]{Ye5} A. Yekutieli,
  Twisted Deformation Quantization of Algebraic Varieties,
  math.AG/0801.3233.
\bibitem[Ye6]{Ye6} A. Yekutieli,
  Central Extensions of Gerbes,
  Advances in Mathematics {\bf 225} (2010), 445-486.
\bibitem[Ye7]{Ye7} A. Yekutieli,
  MC Elements in Pronilpotent DG Lie Algebras
  Eprint arXiv:1103.1035.
\end{thebibliography}

\end{document}